\documentclass[12pt]{article}%
\usepackage{amsfonts}
\usepackage{amssymb}
\usepackage{amsmath}
\usepackage[latin1]{inputenc}
\usepackage{indentfirst}
\usepackage{graphicx}%
\providecommand{\U}[1]{\protect\rule{.1in}{.1in}}
\newtheorem{Theo}{Theorem}
\newtheorem{Lemm}{Lemma}

\newtheorem{Rema}{Remark}

\numberwithin{equation}{section}
\begin{document}

\title{ On the Solution of Locally Lipschitz BSDE Associated to Jump Markov Process.}
\author{K. Abdelhadi, N. Khelfallah\thanks{University of Biskra, Laboratory of Applied
Mathematics, Po.\ Box 145 Biskra (07000), Algeria. E-mail addresses:
abdelhadi\_khaoula@yahoo.com, n.khelfallah@univ-biskra.dz.} }
\date{\today}
\maketitle

\begin{abstract}
In this study, we consider a class of backward SDE driven by jump Markov
process. An existence and uniqueness result to this kind of equations is
obtained in a locally Lipschitz case. We essentially approximate the initial
problem by constructing a convenient sequence of globally Lipschitz BSDEs
having the existence and the uniqueness propriety. Then, we show, by passing
to the limits, the existence and uniqueness of a solution to the initial
problem. After that, a stability theorem is also proved in the local Lipschitz
setting. Applying the aforementioned result, we give an application to
European option pricing with constraint.


{${Keywords}:$} Backward stochastic differential equations, Jump Markov
process, Random measure.

\end{abstract}

\section{Introduction}

The history of backward stochastic differential equations driven by continuous
Brownian motion goes back to the work of J.M. Bismut $\cite{bismut},$ in 1973.
However, the theory that deals with this type of equations was developed in
1990\ by E. Pardoux and S. Peng, in their original paper $\cite{PARDOUX}$. In
$\cite{PARDOUX}$, the authors, announced the first result of existence and
uniqueness concerning the nonlinear case. From this work, many authors attempt
to relax the regularity of the generator with respect\ to the stat variables,
this consists in trading of the Lipschitz condition by imposing more weaker
assumptions. Let us mention that there are a list of large literature in this
respect, see for example
$\cite{Auguste,Kheled,bahlalii,bahlali2,HAMADENE,LEPELTIER,PARDOUX},$ for the
continuous Brownian case; we also refer the reader to\ $\cite{barles,
becherer, confortola1,CONFORTOLA2, GREPEY, KAZI1,KHARROUBI}$ and the
references therein for BSDEs driven by random jumps processes.

Since the aim of this paper is to study locally Lipschitz setting for one kind
of BSDE with jump Markov process, we present some existing results that go in
this direction and are investigated for BSDEs driven by Brownian motion
(without Markov jump part). The first paper treats the\ locally Lipschitz
case, is given by S. Hamadene $\cite{HAMADENE},$ in 1996; in that paper, the
author studied a kind of one dimensional BSDEs with locally Lipschitz
generator satisfying a suitable growing condition assuming that the terminal
condition is bounded. Then, in Bahlali $\cite{Kheled}$, the previous result
has been generalized in the multidimensional case with both local assumption
on the coefficient and only square integrable terminal data. Subsequently, the
same results have been addressed in Bahlali $\cite{bahlalii}$ for BSDEs driven
by Teugel's martingales and an independent Brownian motion. After that, an
existence and uniqueness result to backward stochastic nonlinear Volterra
integral equations under local Lipschitz continuity condition on the drift,
has been investigated in A. Auguste and N. Modeste $\cite{Auguste}.$

More recently, BSDEs driven by a random measure related to a pure jump Markov
process have been studied by F. Confortola, M. Fuhrman in
$\cite{confortola1,CONFORTOLA2}$, where they provided existence and uniqueness
results for such equations. Furthermore, they applied these results to study
nonlinear variants of the Kolmogorov equation of the Markov process and to
solve optimal control problems related to this topic. Herein, we are going to
extend the contribution of the important papers $\cite{Kheled,confortola1}$ by
involving random measure associated to jump Markov process in the state BSDE.

Motivated by the above results, we study a class of backward stochastic
differential equations driven by a random measure associated to jump Markov
process. We tackle an existence and uniqueness result on top of a stability
propriety for the solution of the following type of backward stochastic
differential equation.%
\begin{equation}
Y_{s}=h(X_{T})+\int\limits_{s}^{T}f(r,X_{r},Y_{r},Z_{r}(.))dr-\int
\limits_{s}^{T}\int\limits_{\Gamma}Z_{r}(y)q(dr\text{ }dy),\text{\ \ \ }%
s\in\left[  t,T\right]  . \label{1}%
\end{equation}

Where $q(dt$ $dy)$ is a random measure, $X$ is a normal jump Markov process
and $h(X_{T})$ is the terminal condition. Noting that $Z$, under some
appropriate measurability conditions, is a random field on $\Gamma.$ The
generator $f$ depends on $Y$ and $Z$ in a general functional way. An adapted
solution to BSDE $(\ref{1})$ (if there exists) is a couple $(Y,Z)$ which
belongs to $\mathcal{B}$ and satisfies equation $(\ref{1}).$

In our setting, two main results are established. The first result consists in
proving the existence and uniqueness of the solution to the BSDE$(\ref{1})$
under the locally Lipschitz condition satisfied by the generator $f$. In this
case, we prove our main result assuming that the terminal condition is square
integrable, the generator $f$ satisfies the polynomial growth condition and
the Lipschitz condition in the ball $B(0,M)$ with a Lipschitz constant
inferior or equal to $L+\sqrt{log(M)},$ where $L$ is a universal positive
constant. As a second result, we prove a stability of solutions under the same conditions.

The remainder of the current paper is organized as follows. In Section2, we
give a brief introduction into jump Markov process theory and we recall an
existence and uniqueness Theorem for BSDEs with globally lipschitz
coefficients. In Section 3, we study BSDEs with locally Lipschitz
coefficients. \ To illustrate our theoretical results, an example in finance
is treated in Section 4.

\section{Preliminaries.}

\subsection{Overview of Jump Markov Process}

Throughout this paper, the real positive number $T$ stands for horizon, and
$\left(  \Omega,\mathcal{F},\mathbb{P}\right)  $ stands for a complete
probability space.\ We define $\left(  \Gamma,\mathcal{E}\right)  $ as a
measurable space such that $\mathcal{E}$ contains all one-point sets; we also
define $\Delta$ as a point not included in $\Gamma$. For a given a normal jump
Markov process $X$, we denote by $\mathbb{F}^{t}:=\sigma\left(  X_{r},r\leq
s\right)  \vee\mathcal{N}$ the filtration $\left(  \mathcal{F}_{\left[
t,s\right]  }\right)  _{s\in\left[  t,\infty\right[  }$, generated by the
process $X$ and completed by $\mathcal{N},$ where $\mathcal{N}$ is the
totality of $\mathbb{P}$--null sets. We assume further that the jump Markov
process $X$ satisfies the following conditions:

\begin{enumerate}
\item $\mathbb{P}^{t,x}(X_{t}=x)=1$ for every $t\in\left[  0,\infty\right[  $,
$x\in\Gamma.$

\item For every $0\leq t\leq s$ and $A\in\mathcal{E}$ the function
$x\rightarrow\mathbb{P}^{t,x}(X_{s}\in A)$ is $\mathcal{E}$-measurable.

\item For every $0\leq r\leq t\leq s,$ $A\in\mathcal{E}$ we have
$\mathbb{P}^{r,x}(X_{s}\in A\mid\mathcal{F}_{\left[  r,t\right]  }%
)=\mathbb{P}^{t,X_{t}}(X_{s}\in A).$ $\mathbb{P}^{r,x}$-a.s.

\item All the trajectories of the pure jump process $X$ have right limits when
$\Gamma$ is endowed with its discrete topology (the one where all subsets are
open). In other words, for every $\omega\in\Omega$ and $t\geq0$ there exists
$\delta>0$ such that $X_{s}(\omega)=X_{t}(\omega)$ for $s\in\lbrack
t,t+\delta]$.

\item For every $\omega\in\Omega$ the number of jumps of the trajectory
$t\rightarrow X_{t}(\omega)$ is finite on every bounded interval, which
implies that $X$ is non explosive process.
\end{enumerate}

Let $\mathcal{P}$ be the predictable $\sigma$ -algebra, and \noindent\ $Prog$
be the progressive $\sigma$ -algebra on $\Omega$ $\times\left[  0,\infty
\right[  ,$ the same symbols will also denote the restriction to $\Omega
\times\lbrack t,T].$ We define a transition measure (also called rate measure)
$v(t,x,A),$ $t\in\lbrack0,T],$ $x\in\Gamma,$ $A\in\Gamma$ from $[0,\infty
)\times\Gamma$ to $\Gamma$, such that$\underset{t\in\lbrack0,T],\text{ }%
x\in\Gamma}{\sup}v(t,x,\Gamma)<\infty$ and $v(t,x,\left\{  x\right\}  )=0.$
Let $T_{n}$ be the jump times of $X$, we consider the marked point process
$(T_{n},X_{T_{n}}),$ and the associated random measure $p(dt$ $dy):=\sum
\limits_{n}$ $\delta_{(T_{n},X_{T_{n}})\text{ }}$on $\left(  0,\infty\right)
\times\Gamma,$ where $\delta$\ stands for the Dirac measure. The compensator
(also called the dual predictable projection) $\tilde{p}$ of $p$ is $\tilde
{p}=v(t,X_{t-},dy)dt$, so that $q(dr$ $dy):=p(dr$ $dy)-v(r,X_{r-},dy)dr$ is an
$\mathbb{F}^{t}-$martingale.

We recall the representation theorem of marked point process martingales; it
is one of the important tools to prove the existence and uniqueness of
solution to BSDE. This theorem states that every integrable martingale adapted
to the natural filtration generated by a jump Markov process, can be written
in terms of stochastic integral with respect to this jump Markov process (for
further information in this subject see for example $\cite{CONFORTOLA2,Gihman}%
$)$.$

In the remainder of this study, we will work on the following spaces

\begin{itemize}
\item $\mathcal{L}^{m}(p)$, $m\in\left[  1,\infty\right[  $ the space of
$\mathcal{P}\otimes\mathcal{E}$-measurable real functions $W_{s}(\omega,y)$
defined on $\Omega\times\left[  t,\infty\right[  \times\Gamma$, such that%
\[
\mathbb{E}\int_{t}^{T}\int_{\Gamma}\left\vert W_{s}(y)\right\vert ^{m}%
p^{t}(ds\text{ }dy)=\mathbb{E}\int_{t}^{T}\int_{\Gamma}\left\vert
W_{s}(y)\right\vert ^{m}v(s,X_{s},dy)ds<\infty.
\]

\item $\mathcal{L}_{loc}^{1}(p),$ the space of the real functions $W$ such
that\ $W\mathrm{1\hspace{-0.04in}I}_{\left]  t,\text{ }\tau_{n}\right]  }%
\in\mathcal{L}^{1}(p)$ for some increasing sequence of $\mathbb{F}^{t}%
-$stopping times $\tau_{n}$ diverging to $+\infty$.

\item $L^{2}(\Gamma,\mathcal{E},v(s,x,dy))$ the space of processes
$z:\Gamma\rightarrow%
%TCIMACRO{\U{211d} }%
%BeginExpansion
\mathbb{R}
%EndExpansion
$ such that%
\[
\left\Vert z(.)\right\Vert =(\int\limits_{\Gamma}\left\vert z(y)\right\vert
^{2}v(s,X_{s},dy))^{\frac{1}{2}}<\infty.
\]

\item $\mathcal{B}^{t,x}$ the space of processes $(Y,Z)$ on $\left[
t,T\right]  $ such that%
\[
\left\Vert (Y,Z)\right\Vert _{\mathcal{B}}^{2}=\mathbb{E}\int\limits_{s}%
^{T}\left\vert Y_{r}\right\vert ^{2}dr+\mathbb{E}\int\limits_{s}^{T}\left\Vert
Z_{r}(.)\right\Vert ^{2}dr<\infty.
\]

\end{itemize}

The space $\mathcal{B}$, endowed with this norm, is a Banach space.

\begin{Rema}
A stochastic integral $\int_{t}^{s}\int_{\Gamma}W_{r}(y)q^{t}(dr$ $dy)$ is a
finite variation martingale if $W\in\mathcal{L}^{1}(p^{t}).$
\end{Rema}

\subsection{BSDE with Globally Lipschitz coefficients}

In this subsection we recall a result of existence and uniqueness of solution
to BSDE $(\ref{1})$ in the\ globally Lipschitz framework. Let us first
introduce the following Hypothesis:

\noindent\textbf{Hypothesis}1

\begin{enumerate}
\item[$\mathbf{H}_{1}$] The final condition $h:$ $\Gamma\rightarrow%
%TCIMACRO{\U{211d} }%
%BeginExpansion
\mathbb{R}
%EndExpansion
$ is $\mathcal{F}_{\left[  t,T\right]  }$-measurable function satisfies
$\mathbb{E}\left\vert h(X_{T})\right\vert ^{2}<\infty.$

\item[$\mathbf{H}_{2}$] For every\ $s\in\left[  0,T\right]  ,$ $x\in\Gamma,$
$r\in%
%TCIMACRO{\U{211d} }%
%BeginExpansion
\mathbb{R}
%EndExpansion
,$ $f(s,x,r,.)$ is a mapping \newline$L^{2}(\Gamma,\mathcal{E}%
,v(s,x,dy))\rightarrow%
%TCIMACRO{\U{211d} }%
%BeginExpansion
\mathbb{R}
%EndExpansion
.$

\item[$\mathbf{H}_{3}$] For every bounded and $\mathcal{E}$-measurable
function $z:\Gamma\rightarrow%
%TCIMACRO{\U{211d} }%
%BeginExpansion
\mathbb{R}
%EndExpansion
,$ the mapping $(s,x,r)\rightarrow$ $f(s,x,r,z(.))$ is $B(\left[  0,T\right]
)\otimes\mathcal{E}\otimes B(%
%TCIMACRO{\U{211d} }%
%BeginExpansion
\mathbb{R}
%EndExpansion
)$-measurable.

\item[$\mathbf{H}_{4}$] There exist $L\geq0,$ $\acute{L}\geq0$ such that for
every\ $s\in\left[  0,T\right]  ,$ $x\in\Gamma,$ $r,$ $r^{^{\prime}}\in%
%TCIMACRO{\U{211d} }%
%BeginExpansion
\mathbb{R}
%EndExpansion
,$ $z,$ $z^{^{\prime}}\in L^{2}(\Gamma,\mathcal{E},v(s,x,dy)),$%
\[
\left\vert f(s,x,r,z(.))-f(s,x,\acute{r},\acute{z}(.))\right\vert \leq
\acute{L}(\left\vert r-\acute{r}\right\vert )+L\left\Vert z(.)-\acute
{z}(.)\right\Vert .
\]

\item[$\mathbf{H}_{5}$] $\mathbb{E}\int\limits_{t}^{T}\left\vert
f(s,X_{s},0,0)\right\vert ^{2}ds<\infty.$
\end{enumerate}

Noting that, under the assumptions ($\mathbf{H}_{1}$), ($\mathbf{H}_{2}$) and
($\mathbf{H}_{3}$), Lemma(3.2) in $\cite{CONFORTOLA2},$ shows that the mapping
$\left(  \omega,s,y\right)  \rightarrow f\left(  s,X_{s^{-}}\left(
\omega\right)  ,y,Z_{s}\left(  \omega,.\right)  \right)  $ is
$\mathcal{P\otimes}B(%
%TCIMACRO{\U{211d} }%
%BeginExpansion
\mathbb{R}
%EndExpansion
)$-measurable, if $Z\in\mathcal{L}^{2}(p)$. Furthermore, if $Y$ is
Prog-measurable process then, $\left(  \omega,s\right)  \rightarrow f\left(
s,X_{s^{-}}\left(  \omega\right)  ,Y_{s}\left(  \omega\right)  ,Z_{s}\left(
\omega,.\right)  \right)  $ is Prog-measurable.

The following Theorem has been proved in $\cite{CONFORTOLA2}.$ The authors
showed, in the case where the generator does not depend on the state variables
$Y$ and $Z$, \ that BSDE$(\ref{1})$ admits a unique solution. This done thanks
to the representation Theorem for the $\mathbb{F}-$martingales by means the
stochastic integrals with respect to $q$. Using this last result and a fixed
point argument, they proved the existence and uniqueness of global solution to
BSDE$(\ref{1})$ in the case where the generator depends on state variables $Y$
and $Z$.

\begin{Theo}
\label{theorem1} Suppose that \textbf{Hypothesis}1 holds true. Then BSDE
$(\ref{1})$ has a unique solution $\left(  Y,Z\right)  $ which belongs to
$\mathcal{B}$.
\end{Theo}

\section{BSDE with locally Lipschitz coefficients}

The first purpose of this section is to prove existence and uniqueness of
solutions to BSDE $(\ref{1})$ in the case where the generator $f$ is merely
locally Lipschitz$.$ The second purpose is to prove a stability propriety for
the solutions if they exist. For this end, we need to impose some regularity
assumptions on the coefficients. These conditions are gathered and listed in
following basic assumptions:

\begin{enumerate}
\item[$\mathbf{H}_{2.1.}$] $f$ is continuous in $(y,z)$ for almost all
$\left(  t,\omega\right)  .$

\item[$\mathbf{H}_{2.2}.$] There exists $\lambda>0$ and $\alpha\in\left]
0,1\right[  $ such that%
\[
\left\vert f(s,x,y,z(.))\right\vert \leq\lambda\left[  1+\left\vert
y\right\vert ^{\alpha}+\left\Vert z(.)\right\Vert ^{\alpha}\right]  .
\]

\item[$\mathbf{H}_{2.3}.$] For every $M\in%
%TCIMACRO{\U{2115} }%
%BeginExpansion
\mathbb{N}
%EndExpansion
$, there exists a constant $L_{M}>0$ such that%
\[%
\begin{array}
[c]{l}%
\left\vert f(s,x,y,z(.))-f(s,x,\acute{y},\acute{z}(.))\right\vert \leq
L_{M}\left(  \left\vert y-\acute{y}\right\vert +\left\Vert z(.)-\acute
{z}(.)\right\Vert \right)  ,\\
\mathbb{P}\text{-a.s}.\text{ }\forall t\in\lbrack0,T],
\end{array}
\]
and $\forall$ $y,$ $\acute{y},$ $z,$ $\acute{z}$ such that $\left\vert
y\right\vert <M,$ $\left\vert \acute{y}\right\vert <M,$ $\left\Vert
z(.)\right\Vert <M,$ $\left\Vert \acute{z}(.)\right\Vert <M.$

\item[$\mathbf{H}_{2.4}$] The final condition $h:$ $\Gamma\rightarrow%
%TCIMACRO{\U{211d} }%
%BeginExpansion
\mathbb{R}
%EndExpansion
$ is $\mathcal{F}_{\left[  t,T\right]  }$- measurable function satisfies
$\mathbb{E}\left\vert h(X_{T})\right\vert ^{2}<\infty.$
\end{enumerate}

We assume that $f$ satisfies $(\mathbf{H}_{2.1})$ and $(\mathbf{H}_{2.2})$,
then we define the family of semi--norms $\left(  \Phi_{n}\left(  f\right)
\right)  _{n\in%
%TCIMACRO{\U{2115} }%
%BeginExpansion
\mathbb{N}
%EndExpansion
}$%
\begin{equation}
\Phi_{n}\left(  f\right)  =\left(  \mathbb{E}\int_{0}^{T}\underset{\left\vert
y\right\vert ,\left\Vert z(.)\right\Vert \leq n}{\sup}\left\vert
f(s,X_{s},y,z(.))\right\vert ^{2}ds\right)  ^{\frac{1}{2}}. \label{n}%
\end{equation}

\subsection{A priori estimations and results}

In order to prove the existence and uniqueness result, we start by giving the
following very useful three Lemmas. They involve some priori estimates of
solutions for BSDE $(\ref{1})$ on top of some estimates between two solutions.
These lemmas will be a key tool in proving the essential results of the next
section. For later use, we denote respectively BSDE $(\ref{1}),$ $h\left(
X_{T}\right)  $ by BSDE $(f,h)$ and $\xi.$

\begin{Lemm}
$\label{lemme1}$ Let $(Y,Z)$ be a solution of BSDE $(f,\xi)$.

(\textbf{i}) If $f$ satisfies $(\mathbf{H}_{2.2})$ then there exists a
positive constant $C=C(\lambda,\xi,T),$ which depends on $\lambda^{2}$,
$\mathbb{E}\left\vert \xi\right\vert ^{2}$and $T,$ such that for
every\ $s\in\left[  t,T\right]  $, $\mathbb{E}\left\vert Y_{s}\right\vert
^{2}+\mathbb{E}\int\limits_{s}^{T}\left\Vert Z_{r}(.)\right\Vert ^{2}dr\leq
C.$

(\textbf{ii}) Moreover, if the terminal condition $\xi:=h\left(  X_{T}\right)
$ is bounded$,$ there exists a positive constant $K$, such that $\underset
{t\in\left[  0,T\right]  }{\sup}\left\vert Y_{r}\right\vert ^{2}\leq K.$
\end{Lemm}

\noindent\textbf{Proof of Lemma} $(\ref{lemme1})$

First, we prove (\textbf{i}). Using Itô's formula for semimartingales $($cf.
Theorem 32 in $\cite{Protter})$ to $\left\vert Y_{s}\right\vert ^{2}$ and
integrating on the time interval $\left[  s,T\right]  $ we get%
\begin{equation}
\left\vert Y_{s}\right\vert ^{2}=\left\vert \xi\right\vert ^{2}+2\int
\limits_{s}^{T}Y_{r}f_{r}\left(  r,X_{r},Y_{r},Z_{r}\left(  .\right)  \right)
dr-2\int\limits_{s}^{T}\int\limits_{\Gamma}Y_{r^{-}}Z_{r}\left(  y\right)
q(dr\text{ }dy)-\underset{s\leq r\leq T}{\sum}\left\vert \Delta Y_{r}%
\right\vert ^{2}. \label{i}%
\end{equation}

Noting that the process $\left(  \int\limits_{s}^{.}\int\limits_{\Gamma
}Y_{r^{-}}Z_{r}\left(  y\right)  q(dr\text{ }dy)\right)  _{s\in\left[
s,T\right]  }$ is an $\mathbb{F}^{t}-$ martingale. Due the fact that
$Y_{r^{-}}Z_{r}\left(  y\right)  \in\mathcal{L}^{1}(p),$ one can easily check
that the process $\left(  \int\limits_{s}^{.}\int\limits_{\Gamma}Y_{r^{-}%
}Z_{r}\left(  y\right)  q(dr\text{ }dy)\right)  _{s\in\left[  0,T\right]  }$
is an $\mathbb{F}^{t}-$ martingale. Indeed, from Young's inequality and the
fact that $\underset{t\in\lbrack0,T],\text{ }x\in\Gamma}{\sup}v(t,x,\Gamma
)<\infty$, we get
\begin{align*}
\int\limits_{s}^{T}\int\limits_{\Gamma}\left\vert Y_{r^{-}}\right\vert
\left\vert Z_{r}\left(  y\right)  \right\vert v(r,X_{r},dy)dr  &  \leq\frac
{1}{2}\underset{t,x}{\sup}v(t,x,\Gamma)\int\limits_{0}^{T}\int\limits_{\Gamma
}\left\vert Y_{r}\right\vert ^{2}dr\\
+\frac{1}{2}\int\limits_{0}^{T}\int\limits_{\Gamma}\left\vert Z_{r}\left(
y\right)  \right\vert ^{2}v(r,X_{r},dy)dr  &  <\infty.
\end{align*}

In addition, we can rewrite the last term in the equality $(\ref{i})$ as the
following\textbf{,}%
\begin{align}
\underset{s\leq r\leq T}{\sum}\left\vert \Delta Y_{r}\right\vert ^{2}  &
=\int\limits_{s}^{T}\int\limits_{\Gamma}\left\vert Z_{r}(y)\right\vert
^{2}p(dr\text{ }dy)\label{f}\\
&  =\int\limits_{s}^{T}\int\limits_{\Gamma}\left\vert Z_{r}(y)\right\vert
^{2}q(dr\text{ }dy)+\int\limits_{s}^{T}\int\limits_{\Gamma}\left\vert
Z_{r}\left(  y\right)  \right\vert ^{2}v(r,X_{r},dy)dr.\nonumber
\end{align}

Keeping in mind that the stochastic processes\textbf{ }$q$ is an
$\mathbb{F}^{t}-$ martingale, plugging the equality $(\ref{f})$ into
$(\ref{i})$ and taking the expectation$,$\textbf{ }one can get%
\[%
\begin{array}
[c]{l}%
\mathbb{E}(\left\vert Y_{s}\right\vert ^{2})+\int\limits_{s}^{T}%
\int\limits_{\Gamma}\left\vert Z_{r}\left(  y\right)  \right\vert
^{2}v(r,X_{r},dy)dr=\mathbb{E}\left\vert \xi\right\vert ^{2}+2\mathbb{E}%
\int\limits_{s}^{T}Y_{r}f_{r}(r,X_{r},Y_{r},Z_{r}(.))dr\\
s\in\left[  t,T\right]  .
\end{array}
\]

By invoking $(\mathbf{H}_{2.2})$ and using the inequality $\left\vert
y\right\vert ^{\alpha}\leq1+\left\vert y\right\vert $ for each $\alpha
\in\left[  0,1\right[  ,$ we get%
\begin{equation}
\mathbb{E}(\left\vert Y_{s}\right\vert ^{2})+\mathbb{E}\int\limits_{s}%
^{T}\left\Vert Z_{r}(.)\right\Vert ^{2}dr\leq\mathbb{E}\left\vert
\xi\right\vert ^{2}+6\lambda\mathbb{E}\int\limits_{s}^{T}Y_{r}dr+2\lambda
\mathbb{E}\int\limits_{s}^{T}Y_{r}^{2}dr+2\lambda\mathbb{E}\int\limits_{s}%
^{T}Y_{r^{-}}\left\Vert Z_{r}(.)\right\Vert dr. \label{e}%
\end{equation}
we find, thanks to Young's inequality $2xy\leq\varepsilon x^{2}+\frac{y^{2}%
}{\varepsilon}$ with $\varepsilon=1$
\[
\mathbb{E}(\left\vert Y_{s}\right\vert ^{2})\leq\mathbb{E}\left\vert
\xi\right\vert ^{2}+9T+(1+3\lambda^{2})\mathbb{E}\int\limits_{s}^{T}\left\vert
Y_{r}\right\vert ^{2}dr.
\]

From Gronwall's Lemma, we get
\begin{equation}
\mathbb{E}(\left\vert Y_{s}\right\vert ^{2})\leq\left(  \mathbb{E}\left\vert
\xi\right\vert ^{2}+9T\right)  \exp\left(  (1+3\lambda^{2})T\right)  =C_{1}.
\label{**}%
\end{equation}

We turn back to inequality $(\ref{e})$, it follows that, using Young's
inequality once again with $\varepsilon=2$ and the inequality $(\ref{**})$%
\[
\mathbb{E}\int\limits_{s}^{T}\left\Vert Z_{r}(.)\right\Vert ^{2}dr\leq2\left(
\mathbb{E}\left\vert \xi\right\vert ^{2}+9T\right)  +2T(1+4\lambda^{2}%
)C_{1}=C_{2}.
\]

This proves \textbf{(i). }We proceed now to prove \textbf{(ii). }Firstly, by
invoking and replacing the inequality $(\ref{f})$ into $(\ref{i})$, we get%

\begin{align*}
\left\vert Y_{s}\right\vert ^{2}  &  =\mathbb{E}\left\vert \xi\right\vert
^{2}+2\int\limits_{s}^{T}Y_{r}f_{r}\left(  r,X_{r},Y_{r},Z_{r}\left(
.\right)  \right)  dr-2\int\limits_{s}^{T}\int\limits_{\Gamma}Y_{r^{-}}%
Z_{r}\left(  y\right)  q(dr\text{ }dy)\\
&  -\int\limits_{s}^{T}\int\limits_{\Gamma}\left\vert Z_{r}(y)\right\vert
^{2}q(dr\text{ }dy)-\int\limits_{s}^{T}\int\limits_{\Gamma}\left\vert
Z_{r}\left(  y\right)  \right\vert ^{2}v(r,X_{r},dy)dr.
\end{align*}

Taking the conditional expectation with respect to $\mathcal{F}_{\left[
0,s\right]  }$, using Assumption $(\mathbf{H}_{2.2}),$ the inequality
$\left\vert y\right\vert ^{\alpha}\leq1+\left\vert y\right\vert $ for
$\alpha\in\left[  0,1\right[  $ together with Young's inequality, gives%

\[
\left\vert Y_{s}\right\vert ^{2}\leq C+9T+\left(  2\lambda^{2}+2\lambda
\right)  \int\limits_{s}^{T}\mathbb{E}\left(  \left\vert Y_{r}\right\vert
^{2}\mid\mathcal{F}_{\left[  0,s\right]  }\right)  dr.
\]

For any time $t\leq s,$ using once again $\mathbb{E}\left(  \cdot
\mid\mathcal{F}_{\left[  0,t\right]  }\right)  $ in both sides of the previous
inequality, we get%
\[
\mathbb{E}\left(  \left\vert Y_{s}\right\vert ^{2}\mid\mathcal{F}_{\left[
0,t\right]  }\right)  \leq C+9T+\left(  2\lambda^{2}+2\lambda\right)
\int\limits_{s}^{T}\mathbb{E}\left(  \left\vert Y_{r}\right\vert
^{2}\mathcal{F}_{\left[  0,t\right]  }\right)  dr.
\]

Then, Gronwall's Lemma, yields%
\[
\mathbb{E}\left(  \left\vert Y_{s}\right\vert ^{2}\mid\mathcal{F}_{\left[
0,t\right]  }\right)  \leq\left[  C+9T\right]  \exp\left[  \left(
9T+\lambda^{2}+2\lambda\right)  \left(  T-s\right)  \right]  .
\]

In particular, if $t=s,$ we immediately find%

\[
\left\vert Y_{s}\right\vert ^{2}\leq\left(  C+9T\right)  \exp\left[  \left(
9T+\lambda^{2}+2\lambda\right)  \left(  T-s\right)  \right]  .
\]
Hence%
\[
\underset{s\in\left[  0,T\right]  }{\sup}\left\vert Y_{s}\right\vert ^{2}\leq
K_{1}.
\]

\begin{Lemm}
$\label{lemme2}$ Let $f_{1}$ and $f_{2}$ be two functions satisfy
$(\mathbf{H}_{2.1})$ and $(\mathbf{H}_{2.2})$. Let $(Y^{1},Z^{1})$ [resp.
$(Y^{2},Z^{2})$] be a solution of the BSDE $(f_{1},\xi_{1})$ [resp. BSDE
$(f_{2},\xi_{2})$], where $\xi_{1}$ and\textbf{ }$\xi_{2}$\textbf{ }are two
final conditions satisfy $(\mathbf{H}_{2.4})$. Then for every locally
Lipschitz function $f$ and every $M>1$, the following estimates hold\textbf{:
}%
\begin{align*}
\mathbb{E}(\left\vert Y_{r}^{^{1}}-Y_{r}^{^{2}}\right\vert ^{2})  &
\leq\left[  \mathbb{E}\left\vert \xi_{1}-\xi_{2}\right\vert ^{2}+\Phi_{M}%
^{2}(f-f_{2})+\Phi_{M}^{2}(f_{1}-f)\right.  \hspace{-0.02in}\\
&  \left.  +\frac{C(\xi_{1},\xi_{2},\lambda)}{\left(  1+2L_{M}^{2\text{ }%
}\right)  M^{2(1-\alpha)}}\right]  \exp\left[  (4+4L_{M}^{2\text{ }%
})(T-s)\right]  .
\end{align*}%
\[
\mathbb{E}\int\limits_{s}^{T}\left\Vert Z_{r}^{^{1}}(.)-Z_{r}^{^{2}%
}(.)\right\Vert ^{2}dr\leq C(\xi_{1},\xi_{2},\lambda)\left[  \mathbb{E}%
\left\vert \overset{\_}{\xi}\right\vert ^{2}+\left(  \mathbb{E}\int
\limits_{s}^{T}\left\vert Y_{r}^{^{1}}-Y_{r}^{^{2}}\right\vert ^{2}dr\right)
^{\frac{1}{2}}\right]  .
\]
Where $C(\xi_{1},\xi_{2},\lambda)$ is a constant depends only on $\lambda^{2}%
$, $\mathbb{E}\left\vert \xi_{1}\right\vert ^{2}$ and $\mathbb{E}\left\vert
\xi_{2}\right\vert ^{2}.$
\end{Lemm}

\noindent\textbf{Proof of Lemma} $(\ref{lemme2})$

We set $\bar{Y}=Y^{^{1}}-Y^{^{2}},$ $\overset{\_}{Z}=Z^{^{1}}-Z^{^{2}%
},\overset{\_}{f_{s}}=f_{1}(s,X_{s},Y_{s}^{1},Z_{r}^{1}(.))-f_{2}%
(s,X_{s},Y_{s}^{^{2}},Z_{r}^{2}(.)),$ $\bar{\xi}=\xi^{1}-\xi^{2}.$ By Itô's
formula we have%
\begin{equation}
\mathbb{E}(\left\vert \overset{\_}{Y_{s}}\right\vert ^{2})+\mathbb{E}%
\int\limits_{s}^{T}\left\Vert \overset{\_}{Z_{r}}(.)\right\Vert ^{2}%
dr=\mathbb{E}\left\vert \bar{\xi}\right\vert ^{2}+2\mathbb{E}\int
\limits_{s}^{T}\overset{\_}{Y_{r}}\overset{\_}{f_{r}}dr\text{ \ \ }s\in\left[
t,T\right]  . \label{12}%
\end{equation}
For a given $M>1$, let $L_{M}$ be the Lipschitz constant of $f$ in the ball
$B(0,M),$ we define%
\[
D_{M}:=\left\{  \left(  s,\omega\right)  :\left\vert Y_{s}^{1}\right\vert
^{2}+\left\Vert Z_{s}^{1}(.)\right\Vert ^{2}+\left\vert Y_{s}^{2}\right\vert
^{2}+\left\Vert Z_{s}^{2}(.)\right\Vert ^{2}\geq M^{2}\right\}  ,\text{
}\overset{\_}{D}_{M}:=\Omega\backslash D_{M},
\]

This make it possible to rewrite $(\ref{12})$ as the following%

\begin{equation}%
\begin{array}
[c]{l}%
\mathbb{E}(\left\vert \overset{\_}{Y_{s}}\right\vert ^{2})+\mathbb{E}%
\int\limits_{s}^{T}\left\Vert \overset{\_}{Z_{r}}(.)\right\Vert ^{2}%
dr=\mathbb{E}\left\vert \xi\right\vert ^{2}+2\mathbb{E}\int\limits_{s}%
^{T}\overset{\_}{Y_{r}}\overset{\_}{f_{r}}\mathrm{1\hspace{-0.04in}I}_{D_{M}%
}dr\\
+2\mathbb{E}\int\limits_{s}^{T}\overset{\_}{Y_{r}}\overset{\_}{f_{r}%
}\mathrm{1\hspace{-0.04in}I}_{\overset{\_}{D}_{M}}dr\text{\ \ \ }s\in\left[
t,T\right]  .
\end{array}
\label{m}%
\end{equation}
where $\mathrm{1\hspace{-0.04in}I}_{D^{M}}$ stands for the indicator function
of the set $D.$ We proceed to estimate the last term in the previous equality%

\begin{equation}%
\begin{array}
[c]{l}%
2\mathbb{E}\int\limits_{s}^{T}\overset{\_}{Y_{r}}\overset{\_}{f}%
_{r}\mathrm{1\hspace{-0.04in}I}_{\overset{\_}{D}_{M}}dr=2\mathbb{E}%
\int\limits_{s}^{T}\overset{\_}{Y_{r}[}\left(  f_{1}-f\right)  (r,X_{r}%
,Y_{r}^{1},Z_{r}^{1}(.))]\mathrm{1\hspace{-0.04in}I}_{\overset{\_}{D}_{M}}dr\\
+2\mathbb{E}\int\limits_{s}^{T}\overset{\_}{Y_{r}[}\left(  f-f_{2}\right)
(r,X_{r},Y_{r}^{2},Z_{r}^{2}(.))]\mathrm{1\hspace{-0.04in}I}_{\overset{\_}%
{D}_{M}}dr\\
+2\mathbb{E}\int\limits_{s}^{T}\overset{\_}{Y_{r}[}f(r,X_{r},Y_{r}^{1}%
,Z_{r}^{1}(.))-f(r,X_{r},Y_{r}^{2},Z_{r}^{2}(.))]\mathrm{1\hspace{-0.04in}%
I}_{\overset{\_}{D}_{M}}dr\\
=I_{1}+I_{2}+I_{3}.
\end{array}
\label{13}%
\end{equation}

Then, the inequality $2xy\leq x^{2}+y^{2}$ together with the definition of the
semi-norm $(\ref{n}),$ leads to%
\begin{equation}
I_{1}\leq\mathbb{E}\int\limits_{s}^{T}\left\vert \overset{\_}{Y_{r}%
}\right\vert ^{2}dr+\Phi_{M}^{2}(f_{1}-f), \label{A}%
\end{equation}

and%
\begin{equation}
I_{2}\leq\mathbb{E}\int\limits_{s}^{T}\left\vert \overset{\_}{Y_{r}%
}\right\vert ^{2}dr+\Phi_{M}^{2}(f-f_{2}). \label{A'}%
\end{equation}

Since $f$ is $L_{M}$ -Lipschitz in the ball $B(0,M)$, we obtain%
\[
I_{3}\leq2L_{M}\left[  \mathbb{E}\int\limits_{s}^{T}\left\vert \overset
{\_}{Y_{r}}\right\vert ^{2}dr+\mathbb{E}\int\limits_{s}^{T}\left\vert
\overset{\_}{Y_{r}}\right\vert \left\Vert \overset{\_}{Z_{r}}(.)\right\Vert
dr\right]  .
\]
From Young's inequality $2xy\leq\frac{\gamma^{2}}{2}x^{2}+\frac{2}{\gamma^{2}%
}y^{2},$ we find%
\begin{equation}%
\begin{array}
[c]{l}%
I_{3}\leq\left(  L_{M}^{2}+1\right)  \mathbb{E}\int\limits_{s}^{T}\left\vert
\overset{\_}{Y_{r}}\right\vert ^{2}dr+\frac{\gamma^{2}}{2}\mathbb{E}%
\int\limits_{s}^{T}\left\vert \overset{\_}{Y_{r}}\right\vert ^{2}%
dr+\frac{2L_{M}^{2}}{\gamma^{2}}\mathbb{E}\int\limits_{s}^{T}\left\Vert
\overset{\_}{Z_{r}}(.)\right\Vert ^{2}dr,\\
\leq\left(  L_{M}^{2}+1+\frac{\gamma^{2}}{2}\right)  \mathbb{E}\int
\limits_{s}^{T}\left\vert \overset{\_}{Y_{r}}\right\vert ^{2}dr+\frac
{2L_{M}^{2}}{\gamma^{2}}\mathbb{E}\int\limits_{s}^{T}\left\Vert \overset
{\_}{Z_{r}}(.)\right\Vert ^{2}dr.
\end{array}
\label{B}%
\end{equation}

From the inequalities $(\ref{A}),$ $(\ref{A'})$ and $(\ref{B}),$ one can get%
\begin{equation}%
\begin{array}
[c]{l}%
2\mathbb{E}\int\limits_{s}^{T}\overset{\_}{Y_{r}}\overset{\_}{f}%
_{r}\mathrm{1\hspace{-0.04in}I}_{\overset{\_}{D}_{M}}dr\leq\Phi_{M}%
^{2}(f-f_{2})+\Phi_{M}^{2}(f_{1}-f)\\
+\left(  3+L_{M}^{2}+\frac{\gamma^{2}}{2}\right)  \mathbb{E}\int
\limits_{s}^{T}\left\vert \overset{\_}{Y_{r}}\right\vert ^{2}dr+\frac
{2L_{M}^{2}}{\gamma^{2}}\mathbb{E}\int\limits_{s}^{T}\left\Vert \overset
{\_}{Z_{r}}(.)\right\Vert ^{2}dr.
\end{array}
\label{14}%
\end{equation}

Now, we turn out to estimate the second term in the inequality $(\ref{m}).$
Using Young's inequality $2xy\leq\rho^{2}x^{2}+\frac{y^{2}}{\rho^{2}}$, we get%
\[%
\begin{array}
[c]{l}%
2\mathbb{E}\int\limits_{s}^{T}\overset{\_}{Y_{r}}\overset{\_}{f_{r}%
}\mathrm{1\hspace{-0.04in}I}_{D_{M}}dr\leq\rho^{2}\mathbb{E}\int
\limits_{s}^{T}\left\vert \overset{\_}{Y_{r}}\right\vert ^{2}\mathrm{1\hspace
{-0.04in}I}_{D_{M}}dr+\frac{2}{\rho^{2}}\mathbb{E}\int\limits_{s}%
^{T}\left\vert f_{1}(s,X_{s},Y_{s}^{1},Z_{r}^{1}(.))\right\vert ^{2}%
\mathrm{1\hspace{-0.04in}I}_{D_{M}}dr\\
+\frac{2}{\rho^{2}}\mathbb{E}\int\limits_{s}^{T}\left\vert f_{2}(s,X_{s}%
,Y_{s}^{1},Z_{r}^{1}(.))\right\vert ^{2}\mathrm{1\hspace{-0.04in}I}_{D_{M}}dr,
\end{array}
\]
by using the assumption $(\mathbf{H}_{2.2})$ and the inequality $(a+b+c)^{2}%
\leq3\left(  a^{2}+b^{2}+c^{2}\right)  $, a simple calculation shows that
\[%
\begin{array}
[c]{l}%
2\mathbb{E}\int\limits_{s}^{T}\overset{\_}{Y_{r}}\overset{\_}{f_{r}%
}\mathrm{1\hspace{-0.04in}I}_{D_{M}}dr\leq\rho^{2}\mathbb{E}\int
\limits_{s}^{T}\left\vert \overset{\_}{Y_{r}}\right\vert ^{2}\mathrm{1\hspace
{-0.04in}I}_{D_{M}}dr+\frac{6\lambda^{2}}{\rho^{2}}\mathbb{E}\int
\limits_{s}^{T}\left[  2+\left\vert Y_{r}^{1}\right\vert ^{2\alpha}+\left\Vert
Z_{r}^{1}(.)\right\Vert ^{\alpha}\right. \\
\left.  +\left\vert Y_{r}^{2}\right\vert ^{2\alpha}+\left\Vert Z_{r}%
^{2}(.)\right\Vert ^{\alpha}\right]  \mathrm{1\hspace{-0.04in}I}_{D_{M}}dr.
\end{array}
\]

From Lemma $(\ref{lemme1}),$ Holder's inequality and the fact that%
\[
\mathrm{1\hspace{-0.04in}I}_{D_{M}}\leq M^{-2}\left[  \left\vert Y_{s}%
^{1}\right\vert ^{2}+\left\Vert Z_{s}^{1}(.)\right\Vert ^{2}+\left\vert
Y_{s}^{2}\right\vert ^{2}+\left\Vert Z_{s}^{2}(.)\right\Vert ^{2}\right]  ,
\]

we arrive at
\[%
\begin{array}
[c]{l}%
\mathbb{E}\int\limits_{s}^{T}\left\vert Y_{r}^{1}\right\vert ^{2\alpha
}\mathrm{1\hspace{-0.04in}I}_{D_{M}}dr\leq\left(  \mathbb{E}\int
\limits_{s}^{T}\left\vert Y_{r}^{1}\right\vert ^{2}dr\right)  ^{\alpha}\left(
\mathbb{E}\int\limits_{s}^{T}\mathrm{1\hspace{-0.04in}I}_{D_{M}}dr\right)
^{1-\alpha},\\
\leq\left(  \mathbb{E}\int\limits_{s}^{T}\left\vert Y_{r}^{1}\right\vert
^{2}dr\right)  ^{\alpha}\frac{1}{M^{2(1-\alpha)}}\left[  \mathbb{E}%
\int\limits_{s}^{T}(\left\vert Y_{r}^{1}\right\vert ^{2}+\left\Vert Z_{r}%
^{1}(.)\right\Vert ^{2}+\left\vert Y_{r}^{2}\right\vert ^{2}+\left\Vert
Z_{r}^{2}(.)\right\Vert ^{2})dr\right]  ^{1-\alpha},\\
\leq(CT)^{\alpha}\left[  \left(  C+\acute{C}\right)  \left(  T+1\right)
\right]  ^{1-\alpha}\frac{1}{M^{2(1-\alpha)}}.
\end{array}
\]

Applying the same method to each one of the terms $\mathbb{E}\int
\limits_{s}^{T}\left\vert Y_{r}^{2}\right\vert ^{2\alpha}\mathrm{1\hspace
{-0.04in}I}_{D_{M}}dr,$\newline$\mathbb{E}\int\limits_{s}^{T}\left\Vert
Z_{r}^{1}(.)\right\Vert ^{\alpha}\mathrm{1\hspace{-0.04in}I}_{D_{M}}dr$ and
$\mathbb{E}\int\limits_{s}^{T}\left\Vert Z_{r}^{2}(.)\right\Vert ^{\alpha
}\mathrm{1\hspace{-0.04in}I}_{D_{M}}dr,$ to get
\begin{equation}
2\mathbb{E}\int\limits_{s}^{T}\overset{\_}{Y_{r}}\overset{\_}{f_{r}%
}\mathrm{1\hspace{-0.04in}I}_{D_{M}}dr\leq\rho^{2}\mathbb{E}\int
\limits_{s}^{T}\left\vert \overset{\_}{Y_{r}}\right\vert ^{2}dr+\frac
{C(\xi_{1},\xi_{2},\lambda)}{\rho^{2}M^{2(1-\alpha)}}. \label{B'}%
\end{equation}

If we choose $\gamma^{2}=2L_{M}^{2}$ and $\rho^{2}=\left(  2L_{M}%
^{2}+1\right)  ,$and plugging the inequalities $(\ref{14})$ and $(\ref{B'})$
into $(\ref{m}),$ we find%
\begin{align*}
\mathbb{E}(\left\vert \overset{\_}{Y_{s}}\right\vert ^{2})  &  \leq
\mathbb{E}\left\vert \bar{\xi}\right\vert ^{2}+\Phi_{M}^{2}(f-f_{2})+\Phi
_{M}^{2}(f_{1}-f)+\frac{C(\xi_{1},\xi_{2},\lambda)}{\left(  1+2L_{M}^{2\text{
}}\right)  M^{2(1-\alpha)}}\\
&  +(4+4L_{M}^{2})\mathbb{E}\int\limits_{s}^{T}\left\vert \overset{\_}{Y_{r}%
}\right\vert ^{2}dr.
\end{align*}

We deduce, using Gronwall's Lemma%
\begin{align*}
\mathbb{E}(\left\vert \overset{\_}{Y_{s}}\right\vert ^{2})  &  \leq\left[
\mathbb{E}\left\vert \bar{\xi}\right\vert ^{2}+\Phi_{M}^{2}(f-f_{2})+\Phi
_{M}^{2}(f_{1}-f)+\frac{C(\xi_{1},\xi_{2},\lambda)}{\left(  1+2L_{M}^{2\text{
}}\right)  M^{2(1-\alpha)}}\right] \\
&  \exp\left[  (4+4L_{M}^{2\text{ }})(T-s)\right]  .
\end{align*}

To prove the second inequality we turn back to the equality $(\ref{12}),$ and
we use Schwartz inequality, to obtain
\[
\mathbb{E}\int\limits_{s}^{T}\left\Vert \overset{\_}{Z_{r}}(.)\right\Vert
^{2}dr\leq C(\xi_{1},\xi_{2},\lambda)\left[  \mathbb{E}\left\vert \bar{\xi
}\right\vert ^{2}+\left(  \mathbb{E}\int\limits_{s}^{T}\left\vert \overset
{\_}{Y_{r}}\right\vert ^{2}dr\right)  ^{\frac{1}{2}}\right]  .
\]
this achieve the proof of Lemma 2.

The following Lemma is the main tool in proving our main results; it allows us
to approximate the initial locally Lipschitz problem by constructing a
sequence of globally Lipschitz BSDEs. By taking advantage of the fact that
each BSDE from this sequence has a unique solution, we can prove, by passing
to the limits, that our initial locally Lipschitz BSDE has also a unique
solution under some suitable conditions. Since the proof of this Lemma uses
similar arguments to those goes as that of Lemma $\left(  4.4\right)  $ in
$\cite{bahlalii},$ we omit it in here\textbf{. }

\begin{Lemm}
$\label{lemme3}$\textbf{ }Let $f$ be a function satisfies $(\mathbf{H}_{2.1}%
)$, $(\mathbf{H}_{2.2})$ and $(\mathbf{H}_{2.3})$. Then, there exists a
sequence of functions $f_{n}$ such that,

(\textbf{i}) For each $n$, $f_{n}$ is globally Lipschitz and satisfying
$\mathbf{(H}_{2.2}\mathbf{)}.$

(\textbf{ii}) $\sup_{n}\left\vert f_{n}(s,x,y,z\left(  .\right)  )\right\vert
\leq\left\vert f(s,x,y,z\left(  .\right)  )\right\vert \leq\lambda\left[
1+\left\vert y\right\vert ^{\alpha}+\left\Vert z(.)\right\Vert ^{\alpha
}\right]  .$ $\mathbb{P}$--a.s, $\forall$\ $s\in\left[  t,T\right]  .$

(\textbf{iii}) For every $M$, $\Phi_{M}(fn-f)\rightarrow0$ as $n\rightarrow
\infty$.
\end{Lemm}

Now we are able to state and prove the main results of this paper. It is
important to point out that the usual localization techniques by means of
stopping time do not work when the generator of BSDE $(\ref{1})$ are merely local.

\subsection{The main Theorems}

\subsubsection{Existence and Uniqueness}

\begin{Theo}
\label{theorem2} Suppose that $(\mathbf{H}_{2.1})$, $(\mathbf{H}_{2.2})$,
$(\mathbf{H}_{2.3})$ and $(\mathbf{H}_{2.4})$ hold true. If there exists a
positive constant $L$ such that $L_{M}\leq L+\sqrt{log(M)},$ the BSDE
$(\ref{1})$ has a unique solution $\left(  Y,Z\right)  $ which belongs to
$\mathcal{B}$.
\end{Theo}

\noindent\textbf{Proof }

Suppose that there exists two solutions of the BSDE $(f,\xi)$:$\ \left(
Y^{1},Z^{1}\right)  $ and $\left(  Y^{2},Z^{2}\right)  .$ The proof of the
uniqueness is straight forward of Lemma $(\ref{lemme2})$ applied with
$f_{1}=f_{2}=f,$ $\xi_{1}=\xi_{2}=h\left(  X_{T}\right)  $,\ and Lebesgue's
dominated convergence theorem$.$

To prove the existence, with the sequence of generators $(f_{n})_{n\in%
%TCIMACRO{\U{2115} }%
%BeginExpansion
\mathbb{N}
%EndExpansion
},$ we define a family of approximating BSDEs obtained by replacing the
generator $f$ in BSDE$\left(  1.1\right)  $ by $f_{n}$ defined in Lemma
$(\ref{lemme3})$
\[
Y_{s}^{n}=h\left(  X_{T}\right)  +\int\limits_{s}^{T}f^{n}(r,X_{r},Y_{r}%
^{n},Z_{r}^{n}(.))dr-\int\limits_{s}^{T}\int\limits_{\Gamma}Z_{r}%
^{n}(y)q(dr\text{ }dy),\text{ \ \ }s\in\left[  t,T\right]  .
\]

In view of Theorem $(\ref{theorem1}),$ the above BSDE has a unique solution,
for each integer $n,$ which will be denoted by $(Y^{n},Z^{n})$. Using similar
arguments in the proof of Lemma $(\ref{lemme1}),$ one can easily find
\begin{equation}
\sup_{n}\mathbb{E}\left(  \left\vert Y_{r}^{n}\right\vert ^{2}+\int
\limits_{s}^{T}\left\Vert Z_{r}^{n}(.)\right\Vert ^{2}dr\right)  \leq C.
\label{b}%
\end{equation}

Noting that for each $n\geq N+1,$ $\left\vert f_{n}(s,x,y,z(.))-f_{n}%
(s,x,\acute{y},\acute{z}(.))\right\vert \leq L_{M}\left(  \left\vert
y-\acute{y}\right\vert +\left\Vert z(.)-\acute{z}(.)\right\Vert \right)  .$ We
split the remain of the proof into the following three steps

\noindent\textbf{Step 1}: In this step, let us first assume that $T$ is small
enough such that $T<\frac{(1-\alpha)}{4}.$ Then, we prove that $(Y^{n},Z^{n})$
is a Cauchy sequence in the Banach space $(\mathcal{B}.\left\Vert .\right\Vert
).$ Let us also assume (without loss the generality) that $L=0,$ so that
$L_{M}\leq\sqrt{\log M}.$ We apply Lemma $(\ref{lemme2})$ to $(Y^{n}%
,Z^{n},f_{n},\xi)$, $(Y^{m},Z^{m},f_{m},\xi),$ to obtain%
\begin{align*}
\mathbb{E}\left\vert Y_{r}^{n}-Y_{r}^{m}\right\vert ^{2}  &  \leq\left[
\Phi_{M}^{2}(f_{n}-f)+\Phi_{M}^{2}(f-f_{m})+\frac{C(\xi_{1},\xi_{2},\lambda
)}{\left(  1+2L_{M}^{2}\right)  M^{2(1-\alpha)}}\right] \\
&  \exp\left[  (4+4L_{M}^{2})T\right]  .
\end{align*}
Since $L_{M}\leq\sqrt{\log M},$ we get%
\begin{align*}
\mathbb{E}\left\vert Y_{r}^{n}-Y_{r}^{m}\right\vert ^{2}  &  \leq
K(M,\alpha)\left[  \Phi_{M}^{2}(f_{n}-f)+\Phi_{M}^{2}(f-f_{m})\right. \\
&  \left.  +\frac{C(\xi_{1},\xi_{2},\lambda)}{\left(  2\log M+1\right)
M^{2(1-\alpha)}}\right]  ,
\end{align*}
such that $K(M,\alpha)=\exp\left(  1-\alpha\right)  M^{(1-\alpha)}.$ Passing
to the limits successively on $n,$ $m,$ $M,$ we obtain%
\[
\mathbb{E}\left\vert Y_{r}^{n}-Y_{r}^{m}\right\vert ^{2}\underset{n,\text{
}m,\text{ }M\rightarrow\infty}{\rightarrow}0.
\]
We use $(\ref{b})$ and the Lebesgue's dominated convergence Theorem, to get
\[
\mathbb{E}\int\limits_{s}^{T}\left\vert Y_{r}^{n}-Y_{r}^{m}\right\vert
^{2}dr\underset{n,\text{ }m\rightarrow\infty}{\rightarrow}0.
\]
And then%
\[
\mathbb{E}\int\limits_{s}^{T}\left\Vert Z_{r}^{n}(.)-Z_{r}^{m}(.)\right\Vert
^{2}dr\leq C(\xi_{1},\xi_{2},\lambda)\left[  \left(  \mathbb{E}\int
\limits_{s}^{T}\left\vert Y_{r}^{n}-Y_{r}^{m}\right\vert ^{2}dr\right)
^{\frac{1}{2}}\right]  \underset{n,\text{ }m\rightarrow\infty}{\rightarrow}0.
\]

The two previous limits imply that $(Y^{n},Z^{n})$ is a Cauchy sequence in the
Banach space $(\mathcal{B}.\left\Vert .\right\Vert )$. That is,
\begin{equation}
\exists(Y,Z)\in\mathcal{B}\text{ such that }\underset{n\rightarrow\infty}%
{\lim}\left\Vert (Y^{n},Z^{n})-\left(  Y,Z\right)  \right\Vert _{\mathcal{B}%
}=0. \label{ineq}%
\end{equation}

\textbf{Step2 }In this step, we assume at first that $T$ is an arbitrary large
time duration. Then, we will prove $(Y^{n},Z^{n})$ is a Cauchy sequence in the
Banach space $(\mathcal{B}.\left\Vert .\right\Vert )$ on the time interval
$\left[  0,T\right]  .$ Firstly, let ($[T_{i},T_{i+1}]$)$_{i=0}^{i=k}$ be a
subdivision of $[0,T],$ such that for any $0\leq i\leq k,$ $|T_{i+1}-$
$T_{i}|\leq\delta$, where $\delta$ is a strictly positive number satisfy
$\delta<\frac{(1-\alpha)}{4}$. Now, for $t\in\left[  T_{k-1},T_{k}\right]  ,$
we consider the following BSDE%
\begin{equation}
Y_{s}^{n}=h(X_{T})+\int\limits_{s}^{T_{k}}f(r,X_{r},Y_{r}^{n},Z_{r}%
^{n}(.))dr-\int\limits_{s}^{T_{k}}\int\limits_{\Gamma}Z_{r}^{n}(y)q(dr.dy).
\label{e1}%
\end{equation}

It is obvious from step1 that,the relation $(\ref{ineq})$ remain valid on the
small interval time $\left[  T_{k-1},T_{k}\right]  .$ Next, for $t\in\left[
T_{k-2},T_{k-1}\right]  ,$ we consider the following BSDE%
\begin{equation}
Y_{s}^{n}=Y_{T_{k-1}}^{n}+\int\limits_{s}^{T_{k-1}}f(r,X_{r},Y_{r}%
,Z_{r}(.))dr-\int\limits_{s}^{T_{k-1}}\int\limits_{\Gamma}Z_{r}(y)q(dr.dy).
\label{e2}%
\end{equation}

Due to fact that, $T_{k-1}$ is an element from the previous subinterval
$\left[  T_{k-1},T_{k}\right]  ,$ $Y_{T_{k-1}}^{n}$converges to $Y_{T_{k-1}}$,
and thus, $(Y^{n},Z^{n})$ is a Cauchy sequence and converges to an element in
$\mathcal{B}$, on the small interval $\left[  T_{k-2},T_{k-1}\right]  .$
Repeating this procedure backwardly for $i=k,...,1,$ we obtain the desired
result on the whole time interval $\left[  0,T\right]  .$

\noindent\textbf{Step 3: }In this step, we prove a global existence and
uniqueness solution to BSDE$\left(  1.1\right)  .$ We have, from step 2, the
following result
\[
\exists(Y,Z)\in\mathcal{B}\text{ such that }\underset{n-\rightarrow\infty
}{\lim}\left\Vert (Y^{n},Z^{n})-\left(  Y,Z\right)  \right\Vert _{\mathcal{B}%
}=0.
\]

It remain to prove that $\int\limits_{s}^{T}f_{n}(r,X_{r},Y_{r}^{n},Z_{r}%
^{n}(.))dr$ converges to $\int\limits_{s}^{T}f(r,X_{r},Y_{r},Z_{r}(.))dr$, in
probability. Denoting $\overset{\_}{Y^{n}}=Y^{n}-Y,$ $\overset{\_}{Z^{n}%
}=Z^{n}-Z$ and we set for $M>1$%
\[
A_{M}^{n}:=\left\{  \left(  s,\omega\right)  :\left\vert Y_{s}\right\vert
+\left\Vert Z_{s}(.)\right\Vert +\left\vert Y_{s}^{n}\right\vert +\left\Vert
Z_{s}^{n}(.)\right\Vert \geq M\right\}  ,\overset{\_}{\text{ }A^{n}}%
_{M}:=\Omega\backslash A_{M}^{n}.
\]
We further, denote by $L_{M}$ the Lipschitz constant of $f$ in the ball
$B(0,M)$ and
\[%
\begin{array}
[c]{l}%
\mathbb{E}\int\limits_{s}^{T}\left\vert f_{n}(r,X_{r},Y_{r}^{n},Z_{r}%
^{n}(.))-f(r,X_{r},Y_{r},Z_{r}(.))\right\vert dr\\
=\mathbb{E}\int\limits_{s}^{T}\left\vert \left(  f_{n}-f\right)
(r,X_{r},Y_{r}^{n},Z_{r}^{n}(.))\right\vert \mathrm{1\hspace{-0.04in}I}%
_{\bar{A}_{M}^{n}}dr\\
+\mathbb{E}\int\limits_{s}^{T}\left\vert \left(  f_{n}-f\right)
(r,X_{r},Y_{r}^{n},Z_{r}^{n}(.))\right\vert \mathrm{1\hspace{-0.04in}I}%
_{A_{M}^{n}}dr\\
+\mathbb{E}\int\limits_{s}^{T}\left\vert f(r,X_{r},Y_{r}^{n},Z_{r}%
^{n}(.))-f(r,X_{r},Y_{r},Z_{r}(.))\right\vert \mathrm{1\hspace{-0.04in}%
I}_{\bar{A}_{M}^{n}}dr\\
+\mathbb{E}\int\limits_{s}^{T}\left\vert f(r,X_{r},Y_{r}^{n},Z_{r}%
^{n}(.))-f(r,X_{r},Y_{r},Z_{r}(.))\right\vert \mathrm{1\hspace{-0.04in}%
I}_{A_{M}^{n}}dr.
\end{array}
\]
Then, since $f$ is $L_{M}$ --locally Lipschitz, $(\mathbf{H}_{2.2})$ and
$\left\vert y\right\vert ^{\alpha}\leq1+\left\vert y\right\vert $ for each
$\alpha\in\left[  0,1\right[  $, we obtain%
\[%
\begin{array}
[c]{l}%
\mathbb{E}\int\limits_{s}^{T}\left\vert f_{n}(r,X_{r},Y_{r}^{n},Z_{r}%
^{n}(.))-f(r,X_{r},Y_{r},Z_{r}(.))\right\vert dr\\
\leq\mathbb{E}\int\limits_{0}^{T}\underset{\left\vert y\right\vert ,\left\Vert
z(.)\right\Vert \leq M}{\sup}\left\vert \left(  f_{n}-f\right)  (r,X_{r}%
,y,z(.))\right\vert dr\\
+2\lambda\mathbb{E}\int\limits_{s}^{T}\left[  4+\left\vert Y_{r}%
^{n}\right\vert +\left\Vert Z_{r}^{n}(.)\right\Vert \right]  \mathrm{1\hspace
{-0.04in}I}_{A_{M}^{n}}dr+L_{M}\left(  \mathbb{E}\int\limits_{s}^{T}\left[
\left\vert \overset{\_}{Y_{r}^{n}}\right\vert +\left\Vert \overset{\_}{Z^{n}%
}_{r}\left(  .\right)  \right\Vert \right]  dr\right) \\
+\lambda\mathbb{E}\int\limits_{s}^{T}\left[  6+\left\vert Y_{r}\right\vert
+\left\Vert Z_{r}(.)\right\Vert +\left\vert Y_{r}^{n}\right\vert +\left\Vert
Z_{r}^{n}(.)\right\Vert \right]  \mathrm{1\hspace{-0.04in}I}_{A_{M}^{n}}dr.
\end{array}
\]

Keeping in mind that%
\[
\mathrm{1\hspace{-0.04in}I}_{A_{M}^{n}}\leq M^{-1}\left[  \left\vert
Y_{s}\right\vert +\left\Vert Z_{s}(.)\right\Vert +\left\vert Y_{s}%
^{n}\right\vert +\left\Vert Z_{s}^{n}(.)\right\Vert \right]  ,
\]

then, Schwartz inequality and Lemma $(\ref{lemme1})$ show that there exists a
constant depends on $\xi_{1},$ $\xi_{2},$ $\lambda$ and $T$ such that the
second and the last term in the previous inequality are bounded by
$\frac{C(\xi_{1},\xi_{2},\lambda)}{M^{\frac{1}{2}}}$ and thus,
\[%
\begin{array}
[c]{l}%
\mathbb{E}\int\limits_{s}^{T}\left\vert f_{n}(r,X_{r},Y_{r}^{n},Z_{r}%
^{n}(.))-f(r,X_{r},Y_{r},Z_{r}(.))\right\vert dr\\
\leq\mathbb{E}\int\limits_{0}^{T}\underset{\left\vert y\right\vert ,\left\Vert
z(.)\right\Vert \leq M}{\sup}\left\vert \left(  f_{n}-f\right)  (r,X_{r}%
,y,z(.))\right\vert dr\\
+L_{M}\left(  \mathbb{E}\int\limits_{s}^{T}\left[  \left\vert \overset
{\_}{Y_{r}^{n}}\right\vert +\left\Vert \overset{\_}{Z}_{r}^{n}\left(
.\right)  \right\Vert \right]  dr\right)  +\frac{C(\xi_{1},\xi_{2},\lambda
)}{M^{\frac{1}{2}}}.
\end{array}
\]
Noting that thanks to Lemma $(\ref{lemme3})\left(  \text{\textbf{iii}}\right)
$ and Schwartz inequality, the first term in the previous inequality tends to
$0$ when $n$ goes to infinity. Then, by using Schwartz inequality and step 1,
one can easily check that \newline$\mathbb{E}\int\limits_{s}^{T}\left[
\left\vert \overset{\_}{Y_{r}^{n}}\right\vert +\left\Vert \overset{\_}{Z^{n}%
}_{r}\left(  .\right)  \right\Vert \right]  dr\underset{n\rightarrow\infty
}{\rightarrow}0.$ Finally, due the fact that the constant $C(\xi_{1},\xi
_{2},\lambda)$ is independent of $M,$ the last term goes to $0$ by sending $M$
to infinity. And therefore, $\int\limits_{s}^{T}f_{n}(r,X_{r},Y_{r}^{n}%
,Z_{r}^{n}(.))dr$ converges to$\int\limits_{s}^{T}f(r,X_{r},Y_{r}%
,Z_{r}(.))dr.$ Theorem is proved.

\subsubsection{Stability of the solutions}

Next we will give a stability Theorem for the solution to BDSE $(f,\xi)$ as a
second main result in this paper. Our starting point is to define a sequence
$(f_{n})_{n\in%
%TCIMACRO{\U{2115} }%
%BeginExpansion
\mathbb{N}
%EndExpansion
}$ of Prog--progressively measurable functions, $\left(  \xi_{n}\right)
_{_{n\in%
%TCIMACRO{\U{2115} }%
%BeginExpansion
\mathbb{N}
%EndExpansion
}}$ a sequence of $\mathcal{F}_{\left[  t,T\right]  }$--measurable and square
integrable random variables$.$ For each integer $n$, we suppose that BSDE
$(f_{n},\xi_{n})$ has a (not necessarily unique) solution $(Y^{n},Z^{n}%
)$\textbf{. }Furthermore\textbf{, }We assume also that $(f_{n},\xi_{n})$
satisfies the following assumptions:

\begin{enumerate}
\item[$\mathbf{H}_{3.1}.$] For every $M$, $\Phi_{M}(f_{n}-f)\rightarrow0$ as
$n\rightarrow\infty.$

\item[$\mathbf{H}_{3.2}.$] $\mathbb{E}\left\vert \xi_{n}-\xi\right\vert
^{2}\rightarrow0$ as $n\rightarrow\infty.$

\item[$\mathbf{H}_{3.3}.$] There exist $\lambda>0$ such that:
\end{enumerate}

$\underset{n}{\sup}\left\vert f_{n}(s,x,y,z\left(  .\right)  )\right\vert
\leq\lambda\left[  1+\left\vert y\right\vert ^{\alpha}+\left\Vert
z(.)\right\Vert ^{\alpha}\right]  ,$ $\mathbb{P}$--a.s,$\forall$\ $s\in\left[
t,T\right]  .$

\begin{Theo}
\label{theorem3} $\left(  \text{\textbf{Stability Theorem}}\right)  $ Suppose
that $(f,\xi)$ satisfies Assumptions $(\mathbf{H}_{2.1})$, $(\mathbf{H}%
_{2.2})$, $(\mathbf{H}_{2.3})$ and $(\mathbf{H}_{2.4})$, $(f_{n},\xi_{n})$
satisfies $\left(  \mathbf{H}_{3.1}\right)  $, $\left(  \mathbf{H}%
_{3.2}\right)  $ and $\left(  \mathbf{H}_{3.3}\right)  $. Then we have%
\[
\underset{n\rightarrow\infty}{\lim}\mathbb{E}\int\limits_{s}^{T}(\left\vert
Y_{r}^{n}-Y_{r}\right\vert ^{2})dr+\mathbb{E}\int\limits_{s}^{T}\left\Vert
Z_{r}^{n}(.)-Z_{r}(.)\right\Vert ^{2}dr=0.
\]
\ 
\end{Theo}

\noindent\textbf{Proof }We apply Lemma $(\ref{lemme2})$ to $(Y,Z,f,\xi)$ and
$(Y^{n},Z^{n},f^{n},\xi^{n}),$ we get%
\begin{align*}
\mathbb{E}(\left\vert Y_{r}^{n}-Y_{r}\right\vert ^{2})  &  \leq\left[
\mathbb{E}\left(  \left\vert \xi_{n}-\xi\right\vert ^{2}\right)  +\Phi_{M}%
^{2}(f_{n}-f)+\frac{C(\xi_{1},\xi_{2},\lambda)}{\left(  1+2L_{M}^{2\text{ }%
}\right)  M^{2(1-\alpha)}}\right] \\
&  \exp\left[  (4+4L_{M}^{2\text{ }})(T-s)\right]  .
\end{align*}%
\[
\mathbb{E}\int\limits_{s}^{T}\left\Vert Z_{r}^{n}(.)-Z_{r}(.)\right\Vert
^{2}dr\leq C(\xi_{1},\xi_{2},\lambda)\left[  \mathbb{E}\left(  \left\vert
\xi_{n}-\xi\right\vert ^{2}\right)  +\mathbb{E}(\int\limits_{s}^{T}\left\vert
Y_{r}^{n}-Y_{r}\right\vert ^{2}dr)^{\frac{1}{2}}\right]  .
\]
Passing to the limits, first on $n$ and next on $M,$ and using Lebesgue's
dominated convergence Theorem, we arrive at
\begin{align*}
&  \mathbb{E}(\int\limits_{s}^{T}\left\vert Y_{r}^{n}-Y_{r}\right\vert
^{2}dr)\underset{n,M\rightarrow0}{\rightarrow}0,\\
&  \mathbb{E}\int\limits_{s}^{T}\left\Vert Z_{r}^{n}(.)-Z_{r}(.)\right\Vert
^{2}dr\underset{n\rightarrow0}{\rightarrow}0.
\end{align*}
The Theorem is proved.

\section{Example in finance}

In this section, we adopt the notations of section 2 and we shall give an
example to illustrate our theoretical results. We consider an application to
European option pricing in the constraint case. Noting that the continuos
Brownian case of this model have been treated in $\cite{el karoui}$, under the
globally Lipschitz setting. We impose here to work with a quite general
semimartingale framework assuming that the wealth process is driven by pure
jump Markov process. Roughly speaking, this kind of models comes naturally as
consequence of the lack of continuity in many real world of applications.
Indeed, the empirical distribution of wealth process tend to deviate from
normal distributions. For instance, due to inspected dusters or huge profits,
or even many successive incidents. We start again with a measurable space
$\left(  \Gamma,\mathcal{E}\right)  ,$ tarnation measure $v$ on $\Gamma,$
satisfying $\underset{t\in\lbrack0,T],\text{ }x\in\Gamma}{\sup}v(t,x,\Gamma
)<\infty$ and $v(t,x,\left\{  x\right\}  )=0.$ The jump Markov process $X$ is
also constructed as described in section 2. Our main object is to prove the
existence of feasible strategy to the following general wealth equation%
\[
\left\{
\begin{array}
[c]{ll}%
-dY_{s} & =g\left(  s,X_{s},Y_{s},\pi_{s}\right)  ds-\int_{\Gamma}\pi
_{s}\sigma_{s}q(ds.dy)\\
Y_{T} & =h\left(  X_{T}\right)  \text{ \ \ }s\in\left[  0,T\right]  .
\end{array}
\right.
\]

such that $g$ is a $%
%TCIMACRO{\U{211d} }%
%BeginExpansion
\mathbb{R}
%EndExpansion
$ valued function defined on $\left[  0,T\right]  \times\Gamma\times%
%TCIMACRO{\U{211d} }%
%BeginExpansion
\mathbb{R}
%EndExpansion
\times%
%TCIMACRO{\U{211d} }%
%BeginExpansion
\mathbb{R}
%EndExpansion
,$ $\pi$ is the portfolio process, $\sigma$ is the volatility traded security,
$\xi:=h\left(  X_{T}\right)  \geq0$ is the contingent claim, which is supposed
to be $\mathcal{F}_{\left[  t,T\right]  }$-measurable and square integrable
random variable. Here $Y_{s}$ is the price of the contingent claim $\xi$ at
time $s.$ Noting that this model can be considered as generalization of
Merton's model in two directions: \textbf{(i) }By letting the generator takes
a general form (not necessarily linear). \textbf{(ii) }By involving the Markov
jump part in the wealth process, so that the dynamic contain a discontinuous part.

In the other hand, if we set $Z\left(  .\right)  :=\sigma\pi,$ and we suppose
that $\sigma^{-1}$ is uniformly bounded, the instantaneous wealth process
becomes%
\begin{equation}
\left\{
\begin{array}
[c]{ll}%
-dY_{s} & =f\left(  s,X_{s},Y_{s},Z\left(  .\right)  \right)  ds-\int_{\Gamma
}Z\left(  y\right)  q(ds.dy)\\
Y_{T} & =h\left(  X_{T}\right)  \text{ \ \ \ }s\in\left[  0,T\right]  .
\end{array}
\right.  \label{j}%
\end{equation}

where $g\left(  s,X_{s},Y_{s},\left(  \sigma_{t}\right)  ^{-1}Z\left(
.\right)  \right)  $=$f\left(  s,X_{s},Y_{s},Z\left(  .\right)  \right)  .$

Furthermore, if the generator $g$ satisfies the hypotheses $(\mathbf{H}%
_{2.1})$, $(\mathbf{H}_{2.2})$ and $(\mathbf{H}_{2.3}),$ Theorem
$(\ref{theorem2})$ confirmed to us that BSDE $(\ref{j})$ has a unique solution
$\left(  Y,Z\right)  .$ Therefore, there exists a unique wealth-portfolio
strategy $\left(  Y,\pi\right)  $ which belongs to $\mathcal{B},$ such that
$\pi=\sigma^{-1}Z$.

It remain to prove that the unique existing strategy is feasible (that is, the
constraint of non-negative holds true: $Y_{s}\geq0$ a.s. $t\in\left[
0,T\right]  .)$. To fulfill this, we need to impose the following sufficient
condition, $f\left(  s,X_{s},0,0\right)  $ $\geq0$. Indeed, to benefit from
the comparison Theorem for globally Lipschitz BSDE, we first construct a
sequence of globally Lipschitz generators $(f^{n})_{n\in%
%TCIMACRO{\U{2115} }%
%BeginExpansion
\mathbb{N}
%EndExpansion
},$ by replacing $f$ in BSDE $(\ref{j})$ by $f^{n}$ defined in Lemma
$(\ref{lemme3})$, to obtain the following family of approximating BSDEs
\begin{equation}
Y_{s}^{n}=h\left(  X_{T}\right)  +\int\limits_{s}^{T}f^{n}(r,X_{r},Y_{r}%
^{n},Z_{r}^{n}(.))dr-\int\limits_{s}^{T}\int\limits_{\Gamma}Z_{r}%
^{n}(y)q(dr\text{ }dy),\text{ \ \ }s\in\left[  t,T\right]  , \label{d}%
\end{equation}
In view of Theorem $(\ref{theorem1}),$ BSDE $(\ref{d})$ has a unique solution
$(Y^{n},Z^{n})$, for each integer $n.$ Since $h\left(  X_{T}\right)  \geq0$,
and $f^{n}\geq0$, the comparison Theorem $($cf.Theorem 3.9 in $\cite{bandini}%
),$ leads to $Y^{n}\geq0.$ Using similar arguments in the proof of Theorem
$(\ref{theorem2}),$ one can easily show that $(Y^{n},Z^{n})$ is a Cauchy
sequence then it converges to $\left(  Y,Z\right)  .$ This implies that
$Y\geq0.$ Moreover if the the contingent claim is bounded, lemma
$(\ref{lemme1})$ \textbf{(ii) }shows that the wealth process is also bounded
in the sense that there exists a positive real constant $K$ such that for any
$0\leq t\leq T,$ $0<Y_{t}\leq K.$

\section{Conclusion}

In this work, we have studied a class of backward stochastic differential
equations driven by a random measure associated to jump Markov process.
Motivated by the work of Confortola F, Fuhrman M $\cite{CONFORTOLA2},$ we have
proved an existence and uniqueness result to this kind of equations by
assuming weaker assumptions on the coefficients. More precisely, we have
treated the locally Lipschitz case by using similar techniques developed in
Bahlali $\cite{Kheled}$ with some suitable changes due to the
difference\ between the processes and the spaces. We note that pretty much of
the technical difficulties coming from the pure jump Markov process part are
due to the fact that its quadratic variation, which can be represented as an
integral with respect to the random measure $p(dt$ $dy),$ is not absolutely
continuous with respect to the Lebesgue measure. To overcome these
difficulties, we use the fact that $q(dr$ $dy):=p(dr$ $dy)-v(r,X_{r},dy)dr$ is
an $\mathbb{F}^{t}-$martingale, where $v(r,X_{r},dy)dr$\ represents the dual
predictable projection of $p(dr$ $dy).$

%\begin{thebibliography}{99}
%\small
%\renewcommand{\baselinestretch}{0.3}

\end{document}